\documentclass[12pt,titlepage,twoside]{article}

\usepackage{amsmath, amssymb}
 \usepackage{graphicx}
\usepackage{epstopdf}
 
\newcommand{\rsh}[1]{ \widehat{\pmb{\mathbb R}}^#1}
\newcommand{\rs}[1]{ {\pmb{\mathbb R}}^#1}
\newcommand{\us}[1]{ {\pmb{\mathbb S}}^#1}

\newcommand{\domh}[1]{ \widehat{\bf #1}}

\newtheorem{theorem}{THEOREM }

\newtheorem{corollary}{COROLLARY}

\newtheorem{definition}{DEFINITION}

\begin{document}

\vspace*{2in}

\begin{center}{ { \large D.H. HAMILTON:}  \\

\vspace*{0.5in}

 {\Large MAPPING  THEOREMS}  \\
\vspace*{0.5in}
{\large Department of Mathematics} \\
{\large University of Maryland} \\{\large 2005}  }  
\end{center}

\newpage

\pagestyle{myheadings}
\markboth{\small D.H.Hamilton}{\small  Mapping Theorems}
\setcounter{page}{1}

Riemannian surfaces  are conformally equivalent to the sphere, plane
or disk, in the  case where the  Beltrami differential is zero by the Uniformization Theorem of Poincar\'{e}-Klein-Koebe, i.e. Hilbert's XXII problem (which actually asks for higher dimensions too),  and  more generally when the Beltrami is bounded, by Ahlfors-Bers \cite{key1} who  use $K-$quasiconformal (QC) mappings, i.e. homeomorphisms  mapping small  balls  to ellipsoids of bounded eccentricity $K$.
In higher dimensions Liouville shows that the conformal mappings are essentially trivial \cite{key23} and so one seeks quasi-orthogonal co-ordinates instead. Lavrentieff (1938)  introduced  quasiconformal mappings of space for PDE.  A geometric-analytic theory, initiated by L\"{o}wner (1959), was  sustained  by the school of Gehring and V\"{a}is\"{a}l\"{a} as well as Soviet mathematicians such as Reshetnyak and Zoric, with  notable contributions by Tukia and Rickman (not to mention Ahlfors, Carleson,  Donaldson and Sullivan). In 1965 Gehring and V\"{a}is\"{a}l\"{a} \cite{key9}  formulated what  Ahlfors in his  review called the main  problem: find a  characterisation of  the $K-$quasiconformal images of the unit ball, i.e.  QC-balls.  For  two dimensions this is the Riemann Mapping Theorem,   proved by Koebe in 1907 (since by Ahlfors-Bers \cite{key1} the QC version is equivalent). In their plenary addresses to the  International Congress both Ahlfors(1978)  and Gehring(1986) gave this as the main  open problem of the theory.

\

Our insight comes from reflections\footnote{sense reversing idempotent homeomorphisms of $\rsh{3}$ onto itself}.  Now according to Smith \cite{key28}  any  reflection $ F $ of the sphere $ \rsh{3} $ has a  set $ \bf T $ of fixed points forming a topological sphere with complement being two disjoint domains $ \bf D, \; D' $ (called the  ``complementary domains''). Smith's conjecture  that $ F $ is topologically conjugate to a Euclidean  reflection was disproved by Bing \cite{key4}  by constructing a ``wild reflection'', i.e. the  complementary domains are not simply connected\footnote{``The Smith Conjecture''   program   gave ''yes"  for   diffeomorphisms \cite{key26} }. As Bing's construction of ``wild reflections" was rather indirect he asked \cite{key5} for an explicit example. We constructed  an explicit    example \cite{key14} which although bih\"{o}lder  is not quasiconformal. It was expected (see Heinonen and Semmes \cite{key22}, and communications from Sullivan) that there exist wild $K$-quasiconformal reflections, however in\cite{key16} we prove:

\begin{theorem}{\rm (QC Smith conjecture)} QC reflections are tame.
\end{theorem}

Remarks: Although $ F $ is topologically conjugate to a Euclidean reflection, the conjugation need  not be QC (unlike two dimensions).

\

A fundamental concept in our theory is  renormalization. For example consider
an injection $ F:  \rsh{m} \rightarrow \rsh{n}$ , for $ n \geq m $.
Its family of renormalisations would be $ \widetilde {F}=  N \circ F \circ L $  where $ L $ is any conformal automorphism of $ \rsh{m} $ and  $ N $ chosen so that  say $(0,0,0) $ and   $ (0,0,1) $ are fixed, i.e. $ \widetilde {F}$ is a normalised rescaling. For $ n= m $ this gives a well known characterization of QC mappings $ F $: the family of  renormalizations of a QC map is a  precompact family in the space of  homeomorphisms, i.e. normalized rescalings of quasiconformal mappings have subsequences converging uniformly to another QC map.   For $ n > m $ we get the so called quasisymmetric (QS) mappings which like  QC maps have the property that triangles are roughly preserved. We shall also renormalize sets $ {\bf S} \subset \rsh{3} $, in this case the family is $\widetilde {\bf S}=  L({\bf S}) $  where the conformal mapping fixes a point of $ \bf S $.   For example the famous quasicircles of Ahlfors are Jordan curves $ \Gamma \subset \rsh{2} $ whose family of renormalizations is precompact in the space of Jordan curves (where the Hausdorff metric between sets is used). By Ahlfors these  are parameterized by QS maps of $ \rsh{1} $, i.e.  `` quasisymmetric 1-spheres''.  In general a quasisymmetric m-sphere is a QS embedding of $ \rsh{m} $ in $ \rsh{n} $. Concentrating on  $ \rsh{3}$
we are concerned with $ m=2 $ the so called QS spheres.  In general tan object   with renormalizations being precompact in some metric  will be called uniform. 

\

The fixed set of a $K$-QC reflection is called a quasireflector. By Ahlfors \cite{key2}   the quasireflectors of $ \rsh{2} $ are the QS circles.  We find that in $ \rs{3} $ however there is  a difference between  topological spheres which are uniform as sets and those which have uniform parameterizations.

 \begin{definition}  We say that  { \it  a flat sphere  $ \bf T $ is uniform  if its family of 
``renormalizations'' $ \widetilde{\bf T } $ is precompact in the space of flat spheres with respect to the hausdorff  metric between compact  sets}\footnote{Equivalently, by Bing's criterion for flat spheres, at all scales $ \bf T $ can be squeezed between ``uniform'' polyhedra from the complementary domains}. 
\end{definition}

As the renormalizations  of a QC reflections is a precompact family
from Theorem 1 we deduce that the fixed set is uniform, the converse is

\begin{theorem}{\rm(see \cite{key17})} $ \bf T $ is the fixed set of a QC reflection iff it is a uniform sphere.
\end{theorem}

Remarks: By work of Tukia and V\"{a}is\"{a}l\"{a} this  also solves the problem of characterising the fixed set of bilipschitz reflections (a question first considered by Poincar\'{e} (1898), see Jones\cite{key24}).

\

Another fundamental question was to characterize ``quasispheres'': the image of the unit sphere under a QC mapping $ F :\rsh{3} \rightarrow \rsh{3} $.  The two dimensional analog again is due to Ahlfors who   showed that  ``quasicircle'' was   necessary and sufficient. Of course in $ \rsh{3} $ any quasisphere is a QS sphere, but the converse is false as   QS spheres  can be wild.  Nor is every uniform sphere also QS.  In fact the uniform sphere $ \gamma \times \rs{1} $ for fractal quasicircle $ \gamma $ is not quasisymmetric.
In other words a topological sphere can be uniform as a set 
but without any uniform parametrization.

\begin{definition}: A topological sphere  ${ \bf T} \subset \rs{3}$ is regular if it  a uniform sphere with  quasisymmetric parametrization.
\end{definition}
This provides the characterization required:
\begin{theorem}{\rm(see \cite{key18})} ${ \bf T} \subset \rs{3}$ is the image of $ \us{2}$ under a  QC mapping  of $ \rsh{3}$ iff $ \bf T $ is a   regular sphere. 
\end{theorem}

Remarks:  As a special case we have the  Ahlfors problem of extending a quasisymmetric mapping $ F : \us{2} \rightarrow \us{2} $ to a quasiconformal mapping  of $\rs{3} $, see  L. Ahlfors \cite{key3}, L. Carleson \cite{key8} and (for higher dimensions) Tukia  \cite{key30}.
Indeed  we find various equivalent conditions, e.g.

\begin{corollary} Quasisphere $ \Leftrightarrow $ uniform and  ``L\"{o}wner''.
\end{corollary}

Remarks: ``L\"{o}wner'' (in the terminology of \cite{key6})means that uniform annulli of $ \bf T $ have uniformly bounded 2-capacity.  In fact being a QS sphere is quite delicate.  In \cite{key6}  it is proved that $ \bf T $ is  QS sphere  iff it is `` locally linearly connected'', ``doubling'' and ``L\"{o}wner''. As  our regular spheres already satisfy the first two conditions the corollary is an immediate consequence of Theorem 2.  Essentially all of this means that $ T $ must be ``geometrically uniform''.

\

The QC Riemann Mapping Theorem is a one sided version of this.  Now in  two dimensions the characterization is that $ \bf D $  is simply connected (with non empty boundary) . One difference between two and three dimensions is that any  topological 2-ball  in $\rsh{2} $ is  ``uniformly simply connected '' (USC), i.e. its renormalizations form a precompact family in the space of topological disks (wrt weak convergence).  However  in $ \rs{3} $  a sequence of  renormalizations of a  topological ball could converge  to a torus.   Now a QC-ball is  USC.  However by the previous discussion we expect  extra boundary conditions.     Gehring  \cite{key12} showed the complement $  \rsh{3} - {\bf D}$ is ``linearly locally connected''(LLC).  

\

We must now refer to the prime-ends introduced by Caratheodory
and generalized higher dimensions by Zoric, see  .   These are the 
``ends''  cut off by compacta with diameters converging to zero.
Now any simply connected proper subdomain of $ \rs{2}$ 
has boundary which is a topological circle (in the prime-end metric).
This is not true in $ \rs{2}$ even for topological balls. However 
Zoric showed that  any QC mapping of the unit ball extends to a homeomorphism of the prime-end boundary $ \hat{\partial}{\bf D} $ which is thus a topological sphere.  More generally we find that a USC-LLC  domain  has prime-end boundary $\partial \domh{D} $ parametrised by a  homeomorphism of $ \us{2}$.   Let  $ \cal B $  be the space of prime-end boundary maps of USC-LLC domains using
the natural prime-end metric. 
\begin{definition}  We say that a USC  domain $ \bf D $ is a ``regular ball''  if  it has  parametization $ H: \us{2} \rightarrow  \partial \domh{D} $ whose   renormalizations  are precompact in  $ \cal B $.
\end{definition}

Now  any QC-ball is a regular ball. The converse is the first of our three versions of the QC mapping theorem:

\begin{theorem} Any regular ball is a QC-ball.
\end{theorem}

\

 Previously  V\"{a}is\"{a}l\"{a} \cite{key32} had the best criterion. He  gave necessary and sufficient conditions for  cylindrical domains $ { \bf D = A} \times \rs{1} \subset \rs{3} $ to be the QC image of the ball.   V\"{a}is\"{a}l\"{a} gave several equivalent conditions but  the one we quote    is  that the (prime-end) boundary 
$ \partial \domh{D} $ must be the   image of  $\us{1} \times \rs{1 }$ under parametrization quasisymmetric with respect to the the inner-length:
\[    \nu (X, \; Y ) = \inf \{ dia(\alpha) : X, \; Y \in  \mbox{connected} \; \alpha  \subset
{\bf D}  \} \; , \]  
i.e $\partial \domh{D} $ is $\nu-$quasisymmetric .

\

Indeed a general concept of quasisymmetry is the actual way we proceed to Theorem 4.  We find that USC domains with LLC boundary  is a  Gromov-Uniform-Tree (GUT), see \cite{key19}.  This is analogous to Gromov's theory of Hyperbolic Trees  which arose in the study of discrete groups, for an exposition see \cite{key13}.  Without giving a  formal definition, a GUT is  uniformly approximated from within by uniform trees of uniform polyhedra.  Then  there is a  metric $ \gamma $ on $ \bf D $ which extends to another
ideal boundary, the so called Gromov boundary which turns out to be
equivalent to  the prime-end boundary $ \partial \domh{D} $. The metric is defined by the weights $ c^{-n} $ where $ n $ is the number of steps on disjoint uniform polyhedra from a fixed special point. In fact (analogous to the case of Gromov Hyperbolic Spaces) QC mappings extend to mappings of the Gromov Boundary quasisymmetric in the Gromov metric. Therefore Theorem 4  is actually proved  by\cite{key19}:

\begin{theorem} QC-ball  $\Leftrightarrow$ USC domain bounded by a $\gamma $-QS sphere.
 \end{theorem}
 
Remarks: One way is relatively easy. The converse problem is to extend the QS mapping to a QC mapping of the interiors.

\

 The $\nu-$metric is more explicit than the gromov metric, so  we generalise the V\"{a}is\"{a}l\"{a}  result to the following sufficient condition:
\begin{corollary} Any USC domain bounded by a  
 $\nu$-QS sphere is a QC-ball.
\end{corollary}

Remarks: One could also use the the ordinary (euclidean) QS maps
to get a weaker sufficient condition. 

\ 

Bonk and Kleiner \cite{key7} gave three conditions for a metric space  to be QS equivalent to $ \us{2}$: ``doubling'', LLC and ``Lowner''.  The first two already hold  for USC-LLC domains.  The   ``L\"{o}wner''  criteria is:  {\it uniform annuli of $\partial \domh{D} $ have bounded 2-capacity $ cap $}, i.e.  for any annulus  
\[ { \bf A} =  \{ X : r < \gamma(X,Y) < 2r  \}  \]
 we have $ a < cap({\bf A}) < b $, for absolute constants $a, b $.  This is a thickness condition.  Note that  2-capacity would then be measured in the fairly implicit gromov metric. Actually Bonk and Kliener define 2-capacity in general  metric spaces via approximating circle packings. Likewise we approximate the boundary by PL surfaces
and obtain 2-capacity on the boundary.  However we show that one may use 
the usual capacities (and call the boundary ``L\"{o}wner'')

\

 Our ``three point condition'' is:

\begin{theorem}A domain $ \bf D $ is a QC-ball iff
\begin{enumerate}
\item $\bf D $ is USC
\item $ \bf \rsh{3} - D$ is LLC
\item $\partial { \bf D} $ is L\"{o}wner
\end{enumerate} 
\end{theorem}

Remarks: All three conditions are necessary. For example a  quasireflector satisfies the first two but is not in general QS equivalent to $ \us{2} $. 

\

Sullivan and  Thurston proved \footnote{we thank Sullivan for telling us about this} that a   domain   inscribed
by round balls is a QC-ball.  For example the ``soap bubble'' domain, obtained by first adding to some initial
ball $ \bf B_0 $ some  packing of its boundary by disjoint balls $\bf B_n $ , then doing the same to the $\bf B_n
-B_0$ and continuing outward so that all new balls are disjoint. (An infinite cylinder is another example).
However there are QC-balls which cannot be inscribed by round  balls. Using quasiballs instead we have the concept of quasi-inscribing. Suppose the boundaries of disjoint  K-quasiballs $ \bf \Omega_j $ meet on common faces bounded by curves $ \gamma(i,j)$ where the K-QS maps $H_j: {\bf \partial \Omega_j} \rightarrow \us{2} $ match up, i.e.  $ H_j |_{{ \gamma}(i,j)} = H_i |_{ {\gamma}(i,j)}$. 
\begin{theorem}  Suppose the union $ \bf D $ of the $ \bf \Omega_j $
and the common faces bounded by the $ \gamma(i,j)$ gives a USC domain then $ \bf D $ is a $K'$-QC ball.
\end{theorem}

Example: The so called ``Manhattan '' domains $ {\bf D} $ are formed by adjoining to the half-space $ x_3 < 0 $ the (maybe infinite) cylinders of the form:
\[ {\bf A_j }\times \{0 \leq x_3 < h_j  \leq \infty \} \]
whose bases $ \bf A_j $ are disjoint squares\footnote{actually
any V\"{a}is\"{a}l\"{a} type cylinders with bounded QC constants will do}. Then any Manhattan domain  is a $K$-QC ball, even the one built over a square grid.

\newpage

{\small

\end{document}